\newcount\ite\ite=1\def\0{\global\ite=1\1}
\def\1{\item{\rm(\romannumeral\the\ite)}\advance\ite1\quad}
\def\phi{\varphi}

\documentclass[12pt]{article}
\usepackage{amssymb}

\font\teneufm=eufm10 scaled \magstep1
\font\seveneufm=eufm7 scaled \magstep1
\font\fiveeufm=eufm5  scaled \magstep1
\newfam\eufmfam

\textfont\eufmfam=\teneufm
\scriptfont\eufmfam=\seveneufm
\scriptscriptfont\eufmfam=\fiveeufm

\usepackage{mathrsfs}

\newfam\msbfam
\font\tenmsb=msbm10 scaled \magstep1  \textfont\msbfam=\tenmsb
\font\sevenmsb=msbm7 scaled \magstep1 \scriptfont\msbfam=\sevenmsb
\font\fivemsb=msbm5 scaled \magstep1  \scriptscriptfont\msbfam=\fivemsb

\def\dd#1{\raise1.5pt\hbox{$\,\partial\!$}/\raise-2.5pt\hbox{$\!\partial#1\,$}}

\def\tilde{\widetilde}
\def\hat{\widehat}
\def\5#1{{\mathcal #1}}

\def\RR{{\mathbb R}}
\def\CC{{\mathbb C}}

\def\NN{{\mathbb N}}

\def\PP{{\mathbb P}}

\def\ra{\rightarrow}

\def\GL{\mathop{\rm GL}\nolimits}

\def\Im{\mathop{\rm Im}\nolimits}

\def\End{\mathop{\rm End}\nolimits}

\def\Ann{\mathop{\rm Ann}\nolimits}

\def\Hom{\mathop{\rm Hom}\nolimits}

 \def\HollowBoxx #1#2#3{{\dimen0=#1 \advance\dimen0 by -#2
       \dimen1=#1 \advance\dimen1 by #3
        \vrule height 0pt depth #3 width #2
       \hskip -#3
       \vrule height #1 depth #3 width #3}}
 \def\LeftContraction{\mathord{\kern1.45pt \HollowBoxx{6pt}{3.5pt}{.4pt}}\,}

 \def\HollowBox #1#2#3{{\dimen0=#1 \advance\dimen0 by -#3
       \dimen1=#1 \advance\dimen1 by #3
        \vrule height #1 depth #3 width #3
        \vrule height 0pt depth #3 width #2
        \hskip -#3}}
 \def\RightContraction{\mathord{\, \HollowBox{6pt}{3.1pt}{.4pt}} \kern1.6pt}

\def\qed{{\hfill $\Box$}}
\newtheorem{theorem}{THEOREM}[section]
\newtheorem{corollary}[theorem]{Corollary}

\newtheorem{example}[theorem]{Example}
\newtheorem{remark}[theorem]{Remark}
\newtheorem{proposition}[theorem]{Proposition}

\newtheorem{definition}[theorem]{Definition}

\begin{document}

\begin{center}


\date{23 July 2010}

\end{center}

\begin{center}
{\Large \bf Isolated Hypersurface Singularities\\ 
and Polynomial Realizations\\
\vspace{0.3cm} 
\hspace{0.4cm}of Affine Quadrics}\footnote{{\bf Mathematics Subject Classification:} 32S25, 16A46 }
\medskip\\
\normalsize G. Fels, A. Isaev, W. Kaup, N. Kruzhilin 
\end{center}

\begin{quotation} \small \sl \noindent Let $V$, $\tilde V$ be hypersurface germs in\, $\CC^m$, each having a quasi-homo\-gene\-ous isolated singularity at the origin. We show that the biholomorphic equivalence problem for $V$, $\tilde V$ reduces to the linear equivalence problem for certain polynomials $P$, $\tilde P$ arising from the moduli algebras of $V$, $\tilde V$. The polynomials $P$, $\tilde P$ are completely determined by their quadratic and cubic terms, hence the biholomorphic equivalence problem for $V$, $\tilde V$ in fact reduces to the linear equivalence problem for pairs of quadratic and cubic forms.
\end{quotation}

\markboth{}{}

\setcounter{section}{0}

\section{Introduction}\label{intro}
\setcounter{equation}{0}

Let ${\mathcal O}_m$ be the local algebra of all holomorphic function germs at the origin in $\CC^m$. For every hypersurface germ $V$ at the origin, denote by $I(V)$ the ideal of all elements of ${\mathcal O}_m$ vanishing on $V$. Let $f$ be a generator of $I(V)$, and consider the complex associative commutative algebra ${\mathcal A}(V)$ defined as the quotient of ${\mathcal O}_m$ by the ideal in ${\mathcal O}_m$ generated by $f$ and all its first-order partial derivatives. The algebra ${\mathcal A}(V)$, called the {\it moduli algebra}\, or {\it Tjurina algebra}\, of $V$, is independent of the choice of $f$, and the moduli algebras of biholomorphically equivalent hypersurface germs are isomorphic. Clearly, ${\mathcal A}(V)$ is trivial if and only if $V$ is non-singular. Furthermore, it is well-known that the algebra ${\mathcal A}(V)$ is of finite positive dimension if and only if the germ $V$ has an isolated singularity (see e.g. \cite{GLS}).

By a theorem due to Mather and Yau (see \cite{MY}), two hypersurface germs $V$, $\tilde V$ in $\CC^m$ with isolated singularities are biholomorphically equivalent if their moduli algebras ${\mathcal A}(V)$, ${\mathcal A}(\tilde V)$ are isomorphic. Thus, given the dimension $m$, the moduli algebra ${\mathcal A}(V)$ determines $V$ up to biholomorphism. Moduli algebras, as well as other related associative and Lie algebras, have been extensively studied, and the literature on the moduli spaces of singularities is quite substantial (see e.g. \cite{He1}, \cite{He2} and references therein).

Among all isolated hypersurface singularities, {\it quasi-homogeneous}\, singularities have been of particular interest. Recall that the origin is called a quasi-homogeneous singularity of a hypersurface germ $V$, if for some (hence for any) generator $f$ of $I(V)$ there exist positive integers $p_1,\dots,p_m,q$ such that, modulo a biholomorphic change of coordinates, $f$ is the germ of a polynomial $Q$ satisfying $Q(t^{p_1}z_1,\dots,t^{p_m}z_m)\equiv t^qQ(z_1,\dots,z_m)$ for all $t\in\CC$. By a theorem due to Saito (see \cite{Sa1}), the singularity of $V$ is quasi-homogeneous if and only if $f$ lies in the Jacobian ideal $J(f)$ in ${\mathcal O}_m$, that is, the ideal generated by all first-order partial derivatives of $f$. Hence, for a quasi-homogeneous singularity, ${\mathcal A}(V)$ coincides with the {\it Milnor algebra}\, ${\mathcal O}_m/J(f)$ for any generator $f$ of $I(V)$.

In this paper we present a new criterion for two moduli algebras ${\mathcal A}(V)$, ${\mathcal A}(\tilde V)$ to be isomorphic provided each of $V$, $\tilde V$ has a quasi-homogeneous isolated singularity. In fact, our criterion works for general Gorenstein algebras over $\CC$ of finite dimension greater than 1 (see Theorem \ref{main1}). Recall that a local commutative associative algebra of finite dimension greater than 1 is Gorenstein if and only if the annihilator $\Ann({\mathcal N})$ of its maximal ideal ${\mathcal N}$ is 1-dimensional (see e.g. \cite{B}, \cite{Hu}). In this paper, we call the maximal ideals of such Gorenstein algebras {\it admissible algebras}\, (see Section 2 for the definition), and to every admissible algebra ${\mathcal N}$ and every projection $\pi$ on ${\mathcal N}$ with range $\Ann({\mathcal N})$ we canonically associate a smooth algebraic hypersurface $S_{\pi}\subset{\mathcal N}$, which by way of a special polynomial map is equivalent to an affine quadric. The hypersurface $S_{\pi}$ is essentially independent of the projection $\pi$, more precisely, for two admissible projections the corresponding hypersurfaces differ by a translation on ${\mathcal N}$. We show that two admissible algebras $\5N$, $\tilde{\5N}$ are isomorphic if and only if the hypersurfaces $S_{\pi}$, $S_{\tilde\pi}$ arising from them are affinely equivalent. If at least one of $S_{\pi}$, $S_{\tilde\pi}$ is affinely homogeneous, then these hypersurfaces are affinely equivalent if and only if they are linearly equivalent. We prove that affine homogeneity of $S_{\pi}$ takes place if the algebra $\5N$ admits a grading. 

Next, the hypersurface $S_{\pi}$ is linearly equivalent to the graph of a certain polynomial $P$, which is completely determined by its quadratic and cubic terms and has vanishing constant and linear terms. This polynomial belongs to a special class of polynomials that we call {\it nil-polynomials}\, (see Definition \ref{defnilpol}). We show that $S_{\pi}$, $S_{\tilde\pi}$ are linearly equivalent if and only if the corresponding nil-polynomials $P$, $\tilde P$ are equivalent up to scale by means of a linear transformation (we call such nil-polynomials {\it linearly equivalent}). It then follows that if at least one of the admissible algebras ${\mathcal N}$, $\tilde{\mathcal N}$ admits a grading, then the equivalence problem for them reduces to the linear equivalence problem for the nil-polynomials $P$, $\tilde P$, which in turn reduces to analyzing the quadratic and cubic terms in these polynomials.

If $V$ has a quasi-homogeneous isolated singularity, then the maximal ideal ${\mathcal N}(V)$ of its moduli algebra ${\mathcal A}(V)$ is admissible, provided ${\mathcal N}(V)$ is non-zero. Furthermore, ${\mathcal N}(V)$ admits a grading (in fact, the existence of a grading on ${\mathcal N}(V)$ characterizes quasi-homogeneous singularities -- see \cite{XY}). Therefore the theory developed in Section \ref{section1} can be applied to ${\mathcal N}(V)$. This is done in Section \ref{section2}, where we obtain that the problem of biholomorphic equivalence for quasi-homogeneous isolated singularities of hypersurface germs $V$, $\tilde V$ reduces to the problem of linear equivalence for the corresponding nil-polynomials $P$, $\tilde P$ (see Theorem \ref{main2}). As one can see from Examples \ref{exe8}, \ref{newexample}, in applications it may be useful to exploit higher-order terms in $P$, $\tilde P$ rather than the quadratic and cubic terms alone.

In Section \ref{section2} we show, in particular, how Theorem \ref{main2} works for simple elliptic singularities of type $\tilde E_8$ (see Example \ref{exe8}).  In the article \cite{Ea}, for every singularity of this type a certain cubic homogeneous polynomial has been constructed, with the property that for biholomorphically equivalent singularities the corresponding polynomials are linearly equivalent. Further, the invariant theory for such cubic polynomials was used to distinguish non-equivalent singularities. Interestingly, the polynomial introduced in \cite{Ea} turns out to be part of the corresponding nil-polynomial $P$ arising from our approach. As we will see in Example \ref{exe8}, the nil-polynomial $P$ can be dealt with by an elementary argument that does not require any invariant theory. This simplification is achieved by means of exploiting higher-order terms in $P$. For simple elliptic hypersurface singularities of the remaining two types $\tilde E_6$, $\tilde E_7$, Theorem \ref{main2} leads to considerations similar to those in \cite{Ea}.
\vspace{0.3cm} 

\noindent {\bf Acknowledgements.} We would like to thank M. Eastwood for many inspiring conversations, and A. Neeman and V. Palamodov for useful discussions. Part of this work was done while the fourth author was visiting the Australian National University. The research is supported by the Australian Research Council.

\section{Admissible Algebras and Polynomial\\ Realizations of Certain Affine Quadrics}\label{section1}
\setcounter{equation}{0}

Throughout this section we assume the base field to be $\CC$. Let ${\mathcal N}$ be an associative commutative algebra of finite dimension, and consider the following descending chain of ideals: ${\mathcal N}^1:={\mathcal N}$, ${\mathcal N}^{j+1}:=\langle {\mathcal N}{\mathcal N}^j\rangle$, where $\langle A\rangle$ denotes the linear span of a subset $A\subset{\mathcal N}$. Also, we let ${\mathcal N\,}^0:=\CC\oplus{\mathcal N}$ be the unital extension of ${\mathcal N}$. Then ${\mathcal N}^j{\mathcal N}^m\subset{\mathcal N}^{j+m}$ for all $j,m\ge 0$. Recall that ${\mathcal N}$ is called {\it nilpotent}\, if ${\mathcal N}^{\nu+1}=0$ for some $\nu\ge 0$, and the minimal $\nu$ with this property is called the {\it nil-index}\, of ${\mathcal N}$. Recall further that the {\it annihilator}\, of ${\mathcal N}$ is defined as $\Ann({\mathcal N}):=\{u\in{\mathcal N}:u\,{\mathcal N}=0\}$. Clearly, if ${\mathcal N}$ is nilpotent of nil-index $\nu>0$, then ${\mathcal N}^{\nu}\subset \Ann({\mathcal N})$, hence $\Ann({\mathcal N})$ is non-trivial. We say that a nilpotent algebra ${\mathcal N}$ is {\it admissible}, if $\dim\Ann({\mathcal N})=1$, in which case one has ${\mathcal N}^{\nu}=\Ann({\mathcal N})$. For an admissible algebra ${\mathcal N}$ its unital extension ${\mathcal N\,}^0$ is a finite-dimensional Gorenstein algebra. Since the maximal ideal of any finite-dimensional local algebra is nilpotent by Nakayama's lemma, admissible algebras are exactly the maximal ideals of Gorenstein algebras of finite dimension greater than 1. 

Let ${\mathcal N}$ be an admissible algebra. Fix a
projection $\pi:{\mathcal N}\ra\Ann({\mathcal N})$ with range
$\Ann({\mathcal N})$ (we call such projections {\it admissible}). For
${\mathcal K}:=\ker\pi$ we have ${\mathcal N}=\Ann({\mathcal
  N})\oplus{\mathcal K}$. We extend $\pi$ to ${\mathcal N\,}^0$ by
setting $\pi(1)=0$ and denote the extended projection by the same
symbol. Consider the bilinear $\Ann({\mathcal N})$-valued form
$b_{\pi}$ on ${\mathcal N\,}^0$ defined as
$$
b_{\pi}(u,v):=\pi(uv),\quad u,v\in {\mathcal N\,}^0.
$$
It is well-known that $b_{\pi}$ is non-degenerate on ${\mathcal N\,}^0$ (see e.g. \cite{Ei}, p. 552 or \cite{He2}, implication (i)$\Rightarrow$(iii) in the proof of Lemma 2.2 on p. 11). For completeness of our exposition we include a proof.

\begin{proposition}\label{nondegen} \sl The bilinear form $b_{\pi}$ is non-degenerate.
\end{proposition}

\noindent {\bf Proof:} The radical $\mathcal R:=\{u\in{\mathcal
  N\,}^0:b_{\pi}(u,{\mathcal N\,}^0)=0\}$ of $b_{\pi}$ is orthogonal
to $\CC\subset{\5N\,}^0$ and $\Ann({\mathcal N})$ implying $\mathcal R\subset{\mathcal
  K}$.  For all $j>0$ we have $${\mathcal R}\mathcal
N^j=0\quad\Longrightarrow\quad\mathcal R\mathcal
N^{j-1}\subset\Ann({\mathcal N})\quad\Longrightarrow\quad{\mathcal
  R}\mathcal N^{j-1}=0\,.$$ By recursion down to $j=1$ we conclude
$\mathcal R=0$.\qed \vspace{0.3cm}

Let ${\mathcal N}$ be an admissible algebra of nil-index $\nu$, and let $\PP({\mathcal N\,}^0)$ be the projectivization of ${\mathcal N\,}^0$. Consider the following projective quadric:
$$
Q_{\pi}:=\Bigl\{[u]\in\PP({\mathcal N\,}^0):b_{\pi}(u,u)=0\Bigr\},
$$
where $[u]$ denotes the point of $\PP({\mathcal N\,}^0)$ represented by $u\in{\mathcal N\,}^0$. The inclusion ${\mathcal N}\subset{\mathcal N\,}^0$ induces the inclusion $\PP({\mathcal N})\subset \PP({\mathcal N\,}^0)$, and we think of $\PP({\mathcal N})$ as the hyperplane at infinity in $\PP({\mathcal N\,}^0)$. Also, we identify $1+{\mathcal N}\subset{\mathcal N\,}^0$ with the affine part of $\PP({\mathcal N\,}^0)$, and introduce the corresponding affine quadric $Q\,'_{\!\pi}:=Q_{\pi}\cap (1+{\mathcal N})$. 

Let $\exp: {\mathcal N}\ra 1+{\mathcal N}$ be the exponential map defined as
$$
\displaystyle\exp(u):=1+\sum_{m=1}^{\infty}\frac{u^m}{m!}\,\,.
$$
Notice that in the above sum one has $u^m=0$ for $m>\nu$, thus the exponential map is in fact a polynomial transformation. It is bijective with polynomial inverse given by
$$
\log(1+u):=\sum_{m=1}^{\infty}\frac{(-1)^{m+1}}{m}u^m,
$$
for $u\in{\mathcal N}$. As shorthand we use the following
notation: for every element $u\in\5N$ and all integers $j,m\ge1$ put
\begin{equation}
u^{(m)}:={u^{m}\over m!}\quad\hbox{and}\quad\exp_{j}(u):=\sum_{m=j}^{\infty}u^{(m)}\;.\label{shorthand}
\end{equation}

Following \cite{FK2} (see also \cite{FK1}), we now consider
$S_{\pi}:=\log(Q\,'_{\!\pi})$. Observe that $S_{\pi}$ is given by the equation
$\pi(\exp_1(2u))=0$, hence it is a smooth
algebraic hypersurface in ${\mathcal N}$ passing through the origin. The hypersurfaces
$S_{\pi}$ depend on the choice of the admissible projection $\pi$. We
describe this dependence in the following proposition.

\begin{proposition}\label{depdendence} \sl Let ${\mathcal N}$ be an
  admissible algebra with admissible projections $\pi$,
  $\tilde\pi$. Then there exists $a\in{\mathcal N}$ with
  $S_{\tilde\pi}=S_{\pi}+a$.
\end{proposition}

\noindent {\bf Proof:} For $\5A:=\Ann(\5N)$ there exists
$\lambda\in\Hom(\5N,\5A)$ with $\lambda(\5A)=0$ and
$\tilde\pi(u)=\pi(u)+\lambda(u)$ for all $u\in\5N$, where $\Hom$ denotes the space of linear operators.  For every
$v\in\5N$ define the multiplication operator
$L_{v}\in\End(\5N)$ by $u\mapsto vu$. It is easy to see that $v\mapsto
\pi\circ L_{v}$ defines a linear isomorphism
$\5N/\5A\cong\Hom(\5N/\5A,\5A)$. This implies $\lambda=\pi\circ L_{c}$
for some $c\in\5N$ with $\pi(c)=0$. But $\pi(c)=0$ forces
$$
\tilde\pi(\exp_1(2u))=\pi((1+c)\exp_1(2u))=
\pi(\exp_1(2(u-a)))
$$ 
for $a:=-1/2\log(1+c)$, that is,
$S_{\tilde\pi}=S_{\pi}+a$ as claimed.\qed

\medskip Next, we show that the affine equivalence class of the
hypersurface $S_{\pi}$ determines the algebra ${\mathcal N}$ up to
isomorphism.

\begin{proposition}\label{prop1prime}\sl Let ${\mathcal N}$, $\tilde {\mathcal N}$ be admissible algebras with admissible projections $\pi$, $\tilde\pi$, respectively. Then ${\mathcal N}$, $\tilde {\mathcal N}$ are isomorphic if and only if the hypersurfaces $S_{\pi}$, $S_{\tilde\pi}$ are affinely equivalent.
\end{proposition}

\noindent{\bf Proof:} The necessity follows from Proposition \ref{depdendence}, so we assume that $S:=S_{\pi}$ and $\tilde S:=S_{\tilde\pi}$ are affinely equivalent. We write the common dimension of ${\mathcal N}$, $\tilde{\mathcal N}$ as $n+1$ and only consider the non-trivial case $n\ge 1$.

Forgetting the complex structure on ${\mathcal N}$, we denote by ${\mathcal N}\,^{\CC}:={\mathcal N}\oplus i{\mathcal N}$ the formal complexification of ${\mathcal N}$. Then ${\mathcal N}\,^{\CC}$ is a nilpotent complex algebra with annihilator of complex dimension two and unital extension $({\mathcal N}\,^{\CC})^{0}\simeq(\RR\oplus\5N)^{\CC}$. By $w\mapsto\bar w:=u-iv$ for all $w=u+iv$ with $u,v\in{\mathcal N}$ we denote the conjugation defining the real form ${\mathcal N}$ of ${\mathcal N}\,^{\CC}$ and extend it to the conjugation defining the real form $\RR\oplus\5N$ of $({\mathcal N}\,^{\CC})^{0}$. Further, we extend $\pi$ to a complex-linear projection on $({\mathcal N}\,^{\CC})^{0}$ whose kernel contains 1. The extended conjugation and projection will be denoted by the same respective symbols. Then 
$
h(w,w'):=\pi(w\bar w')
$
is an annihilator-valued hermitian form on $({\mathcal N}\,^{\CC})^{0}$, which coincides with $b_{\pi}$ on $\RR\oplus{\mathcal N}$.

Consider the following real quadric of codimension 2 in $\PP(({\mathcal N}\,^{\CC})^{0})$:
$$
{\mathcal Q}:=\left\{[w]\in\PP(({\mathcal N}\,^{\CC})^{0}):h(w,w)=0\right\},
$$
where $[w]$ denotes the point of $\PP(({\mathcal N}\,^{\CC})^{0})$ represented by $w\in({\mathcal N}\,^{\CC})^{0}$. It is straightforward to check that the exponential map associated to ${\mathcal N}\,^{\CC}$ transforms the real codimension two submanifold $S+i{\mathcal N}$ of ${\mathcal N}\,^{\CC}$ into the affine part ${\mathcal Q'}:={\mathcal Q}\cap(1+{\mathcal N}\,^{\CC})$ of the quadric ${\mathcal Q}$. Observe that ${\mathcal Q'}$ is a real-analytic Levi non-degenerate minimal CR-submanifold of $1+{\mathcal N}\,^{\CC}$ of real codimension two. In fact, ${\mathcal Q'}$ is linearly equivalent to the real codimension two quadric ${\mathcal Q}_0$ in $\CC^{2n+2}$ given by the equations:
$$
\begin{array}{l}
\displaystyle\Im w_1=\sum_{j=1}^n(|z_j|^2-|z_{j+n}|^2)\,,\\
\vspace{-0.3cm}\\
\displaystyle\Im w_2=\sum_{j=1}^n(z_j\overline{z}_{j+n}+z_{j+n}\overline{z}_j)\,,
\end{array}
$$
where coordinates in $\CC^{2n+2}$ are denoted by $w_1$, $w_2$, $z_j$ for $j=1,\dots,2n$. Indeed, one can see that ${\mathcal Q'}$ and ${\mathcal Q}_0$ are linearly equivalent by choosing coordinates in the complex algebra ${\5N\,}^0$ such that the restriction of the bilinear form $b_{\pi}$ to $\ker\pi\subset\5N$ is given by the identity matrix. Further, we do the same for $\tilde{\5N}$ and obtain a hermitian form $\tilde h$ on $({\mathcal N}\,^{\CC})^{0}$, as well as a real quadric of codimension two $\tilde{\5Q}$ in $\PP(({\mathcal N}\,^{\CC})^{0})$ and the corresponding affine quadric $\widetilde{\5Q'}$. 

Let $A:{\5N}\ra\tilde{\5N}$ be a complex affine map that establishes
equivalence between $S$ and $\tilde S$. We treat $A$ as a real affine
map and extend it to a complex affine map
$A^{\CC}:{\5N}\,^{\CC}\ra\tilde{\5N}\,^{\CC}$. Consider the biholomorphic
transformation $\Phi:=\tilde\exp\circ\, A^{\CC}\circ
\log:1+{\5N}\,^{\CC}\ra 1+\tilde{\5N}\,^{\CC}$, where $\log:=\exp^{-1}$, and $\exp$, $\tilde\exp$ are the exponential maps associated to ${\mathcal N}\,^{\CC}$, $\tilde{\mathcal N}\,^{\CC}$, respectively. Clearly, $\Phi$ is a
polynomial map, has a polynomial inverse, and transforms ${\5Q}\,'$
into $\tilde {\5Q}\,'$. Since both ${\5Q}\,'$ and $\tilde {\5Q}\,'$
are linearly equivalent to the quadric ${\5Q}_0$, the map $\Phi$
induces a polynomial automorphism $F$ of $\CC^{2n+2}$, which has a
polynomial inverse and preserves ${\5Q}_0$. In particular, $F$ is a global CR-automorphism of the quadric ${\5Q}_0$ and hence is
affine (see the elliptic case on pp. 37--38 in \cite{ES}). Therefore $\Phi$ is affine as well. But for $u\in{\mathcal N}$ we have
$$
\Phi(1+u)=\tilde\exp(a)\left(1+L(u)+{1\over2}\Big(L(u)^{2}-L(u^{2})\Big)+\hbox{higher-order terms}\right)\,,
$$
where $a:=A(0)$ and $L$ is the linear part of $A$. This implies that $L(u)^{2}=L(u^{2})$ for all $u\in\5N$, that is, $L:\5N\to\tilde{\5N}$ is an isomorphism of algebras. \qed 
\vspace{0.3cm}

By Proposition \ref{prop1prime}, the question whether two admissible
algebras ${\5N}$, $\tilde {\5N}$ are isomorphic
reduces to the question whether the corresponding hypersurfaces $S_{\pi}$, $S_{\tilde\pi}$ are affinely equivalent. If at least one of $S_{\pi}$, $S_{\tilde\pi}$
is affinely homogeneous, then these hypersurfaces are affinely equivalent if and only if they are linearly equivalent. Indeed, if, for instance, $S_{\pi}$ is affinely homogeneous and $A:{\mathcal N}\ra\tilde{\mathcal N}$ is an affine map that establishes equivalence between $S_{\pi}$, $S_{\tilde\pi}$, then the linear map $A\circ A'$ also establishes equivalence between $S_{\pi}$, $S_{\tilde\pi}$, where $A'$ is an affine automorphism of $S_{\pi}$ such that $A'(0)=A^{-1}(0)$. 

Below we obtain affine homogeneity of every $S_{\pi}$ for a certain class of admissible algebras. This class is fairly large
and contains all algebras occurring in the next section.

Let ${\mathcal N}$ be an admissible algebra with a grading, that is,
$$ 
{\mathcal N}=\bigoplus_{j>0}{\mathcal N}_{j},\quad {\mathcal N}_{j}{\mathcal N}_{m}\subset{\mathcal N}_{j+m}\;,
$$
where ${\mathcal N}_{j}$ are linear subspaces of ${\mathcal N}$. Then $\Ann({\mathcal N})=\5N_d$ for $d:=\max\{j:{\mathcal N}_{j}\ne\{0\}\}$. We denote by
$\pi:{\mathcal N}\ra {\mathcal N}_{d}$ the projection
with kernel ${\mathcal K}:=\bigoplus_{j<d}{\mathcal N}_{j}$. Consider the
polynomial map $f:=\pi\circ\exp_1:{\mathcal N}\ra{\mathcal
  N}_{d}$. Then the submanifold $S:=f^{-1}(0)$ of ${\mathcal N}$ is
the graph of a polynomial map from ${\mathcal K}$ to ${\mathcal
  N}_{d}$. Clearly, $S$ is linearly equivalent to $S_{\pi}$.

The map $f$ can be written as follows. Every point $u\in{\mathcal N}$ has a unique representation $u=u_1+\dots+u_{d}$ with $u_j\in{\mathcal N}_{j}$. Then we have
$$
f(u)=\sum_{||\mu||=d}u^{(\mu)},
$$
where $\mu=(\mu_1,\dots,\mu_d)\in\NN^d$ is a multi-index, $u^{(\mu)}:=u_1^{(\mu_1)}\cdot\dots\cdot u_d^{(\mu_d)}$ and $||\mu||:=\sum_{m=1}^d m\mu_m$. Further, we let $A(f)$ be the group of all affine transformations $g$ of ${\mathcal N}$ such that $f\circ g=f$.

With this notation, we have the following proposition.

\begin{proposition}\label{graded}\sl The group $A(f)$ acts transitively on $S$. 
\end{proposition}

\noindent {\bf Proof:} For every $1\le j<d$ and $\alpha\in{\mathcal N}_{j}$ consider the following vector field on ${\mathcal N}$:
$$
\xi_{\alpha}(u):=\left(d-j-\sum_{m=1}^{d-j}m u_m\right)\alpha.
$$
Let $L_{\xi_{\alpha}}$ be the linear differential operator corresponding to this vector field. Applying $L_{\xi_{\alpha}}$ to the map $f$ and setting $v^{(-1)}:=0$ for any $v\in{\mathcal N}$, we obtain  
$$
\begin{array}{l}
\displaystyle L_{\xi_{\alpha}}f=\sum_{||\mu||=d}\Biggl[(d-j)\alpha u_1^{(\mu_1)}\cdot\dots\cdot u_j^{(\mu_j-1)}\cdot\dots\cdot u_d^{(\mu_d)}\;-\\
\vspace{-0.3cm}\\
\hspace{2.5cm}\displaystyle\sum_{m=1}^{d-j}m(\mu_{m}+1)\alpha u_1^{(\mu_1)}\cdot\dots\cdot u_m^{(\mu_m+1)}\cdot\dots\cdot 
u_{j+m}^{(\mu_{j+m}-1)}\cdot\dots\cdot u_d^{(\mu_d)}\Biggr]\\
\vspace{-0.3cm}\\
\hspace{1cm}=\;\displaystyle\sum_{\nu}\left(d-j-\sum_{m=1}^{d-j}m\nu_m\right)\alpha u^{(\nu)}\;,
\end{array}
$$
where the last sum is taken over some set ${\mathcal S}$ of multi-indices $\nu=(\nu_1,\dots,\nu_d)\in\NN^d$. It is easy to see that every $\nu\in{\mathcal S}$ satisfies $\sum_{m=1}^dm\nu_m=d-j$. In particular, we have $\nu_m=0$ for $m>d-j$, hence $L_{\xi_{\alpha}}f=0$. 

Let ${\mathfrak a}$ be the linear span of all vector fields $\xi_{\alpha}$ on ${\mathcal N}$, where $\alpha\in{\mathcal N}_{j}$ and $1\le j<d$. It is easy to see that ${\mathfrak a}$ is a Lie algebra of dimension $\dim{\mathcal K}=\dim S$ and the evaluation map $\varepsilon_u:{\mathfrak a}\ra T_u\,{\mathcal N}={\mathcal N}$, $\xi\mapsto\xi_u$ is injective for every $u\in{\mathcal N}$. Since $L_{\xi}f=0$ for all $\xi\in{\mathfrak a}$ and $S$ is a closed complex submanifold of ${\mathcal N}$, the group $A(f)$ acts on $S$ transitively, as required.\qed
\vspace{0.3cm}\\

Propositions \ref{depdendence} and \ref{graded} imply that for every admissible algebra $\5N$ and every admissible projection $\pi$ on $\5N$ the hypersurface $S_{\pi}$ is affinely homogeneous if $\5N$ admits a grading.

\begin{corollary}\label{lineqalg} \sl Let ${\mathcal N}$, $\tilde {\mathcal N}$ be admissible algebras with admissible projections $\pi$, $\tilde \pi$, respectively. Assume that at least one of $\5N$, $\tilde{\5N}$ admits a grading. Then the algebras ${\mathcal N}$, $\tilde {\mathcal N}$ are isomorphic if and only if the hypersurfaces $S_{\pi}$, $S_{\tilde\pi}$ are linearly equivalent.
\end{corollary}

\begin{remark}\label{nongraded} \rm Not every admissible algebra admits a grading. An example is given by the maximal ideal of the algebra $A$ defined in Example 1.4 of \cite{CK} (see also Remark 3.3 therein). A slightly varied example with the same property is the maximal ideal $\5N$ of $\5O_2/(z_1^{3}z_2,z_1^{5},z_1z_2^{3}+z_1^{3},z_1^{2}z_2^{2}+z_2^{4})$.  With computer aid it can be seen that for every admissible projection $\pi$ on $\5N$ the corresponding hypersurface $S_{\pi}$ is affinely homogeneous. Since we do not know of any admissible algebra that fails to have this property, it remains open whether the existence of a grading is a superfluous  condition for affine homogeneity of $S_{\pi}$.
\end{remark}

So far, we have reduced the isomorphism problem for given admissible
algebras $\5N$, $\tilde{\5N}$ to the affine equivalence problem for the
associated hypersurfaces $S_{\pi}$, $S_{\tilde\pi}$, and even to the linear equivalence problem for these hypersurfaces, if at least one of $\5N$, $\tilde{\5N}$ admits a grading.  Graded admissible
algebras will play a prominent role in applications in the next
section. In the following we give a further reduction of the isomorphism
problem for admissible algebras to the equivalence problem for certain polynomials.

For every complex vector space $W$ of finite dimension we denote by
$\CC[W]$ the algebra of all $\CC$-valued polynomials on $W$.

\begin{definition}\label{defnilpol} \rm A polynomial $P\in\CC[W]$ is called a {\it
      nil-polynomial}\, if there exists an admissible algebra ${\mathcal
      N}$, a linear form $\omega:{\mathcal N}\ra\CC$ and a linear
    isomorphism $\phi:W\ra\ker\omega$ such that
    $\omega(\Ann({\mathcal N}))=\CC$ and
    $P=\omega\circ\exp_{2}\circ\phi$. Two nil-polynomials
    $P\in\CC[W]$, $\tilde P\in\CC[\tilde W]$ are called {\it linearly equivalent} if there exists a linear isomorphism $g:W\to\tilde W$ and
    $r\in\CC^{*}$ such that $P=r\cdot\tilde
    P\circ g$.
\end{definition}

If $\omega$ is a linear form as in Definition \ref{defnilpol}, there exists a unique admissible projection $\pi$ on $\5N$ with
$\ker\pi=\ker\omega$. Conversely, for every admissible projection $\pi$ on $\5N$ there exist a unique, up to a scalar factor, linear form $\omega$ with this property. If the kernels of $\omega$ and $\pi$ coincide, then for every nil-polynomial $P$ related to $\omega$ as in Definition \ref{defnilpol}, the hypersurface $S_{\pi}$ is linearly equivalent to the graph of $P$ in $\CC\times W$. 

The nil-polynomial $P$ has a unique decomposition
$$
\displaystyle P=\sum_{\ell=2}^{\nu}P^{[\ell]}\;,\qquad P^{[\ell]}(x)=\omega(\phi(x)^{(\ell)})
$$
(see (\ref{shorthand})), where every $P^{[\ell]}\in\CC[W]$ is
homogeneous of degree $\ell$ and $\nu$ is the nil-index of $\5N$. The
quadratic form $P^{[2]}$ is non-degenerate on $W$, and $P^{[\nu]}\ne0$,
provided $W\ne\{0\}$, that is, $\dim\5N\ge2$. For every $\ell\ge2$ there
exists a unique symmetric $\ell$-linear
form 
$$
\omega^{}_{\ell}:W^{\ell}\to\CC\quad\hbox{with}\quad\omega^{}_{\ell}(x,\dots,x)=\ell!\,P^{[\ell]}(x)
$$
for all $x\in W$. Clearly
\begin{equation}
\omega_{\ell}(x^1,\dots,x^{\ell})=\omega\Bigl(\varphi(x^1)\cdots\varphi(x^{\ell})\Bigr)\label{expressomega}
\end{equation}
for all $x^1,\dots,x^{\ell}\in W$. Using $\omega_{2}$ and
$\omega_{3}$ we define a commutative product $(x,y)\mapsto x\cdot y$
on $W$ by requiring the identity
$$
\omega_{2}(x\cdot y,z)=\omega_{3}(x,y,z)
$$
to hold for all $x,y,z\in W$.

\begin{proposition}\label{determ}\sl For $\ell\ge 2$ and all $x^1,\dots,x^{{\ell}+1}\in W$  we have
$$
\omega_{\ell+1}(x^1,\dots,x^{{\ell}+1})=\omega_{\ell}(x^1\cdot x^2,x^3,\dots,x^{\ell+1})\;.
$$
\end{proposition}

\noindent {\bf Proof:} By (\ref{expressomega}) it is sufficient to
show that $\varphi(x)\varphi(y)-\varphi(x\cdot y)\in\Ann({\mathcal N})$ for all $x,y\in W$. For $\pi$ related to $\omega$ as above, we have
$$
\begin{array}{l}
\omega\Bigl(\varphi(x\cdot y)\varphi(z)\Bigr)=\omega_{2}(x\cdot y,z)=\omega_{3}(x,y,z)=\omega\Bigl(\varphi(x)\varphi(y)\varphi(z)\Bigr)=\\
\vspace{-0.3cm}\\
\omega\left(\left[\varphi(x)\varphi(y)-\pi\Bigl(\varphi(x)\varphi(y)\Bigr)\right]\varphi(z)\right)
\end{array}
$$ 
for all $x,y,z\in W$. Put
$$
a:=\varphi(x\cdot y)-\left[\varphi(x)\varphi(y)-\pi\Bigl(\varphi(x)\varphi(y)\Bigr)\right].
$$
Since $\omega(a)=0$ and $b_{\pi}(a,\5N)=0$, arguing as in the proof of Proposition \ref{nondegen} we obtain that $a=0$. Hence
$\varphi(x)\varphi(y)-\varphi(x\cdot
y)=\pi\Bigl(\varphi(x)\varphi(y)\Bigr)\in\Ann({\mathcal N})$, as
required.\qed \vspace{0.3cm}

The recursion formula in Proposition \ref{determ} allows one to recover every $P^{[\ell]}$ from $P^{[2]}$ and
$P^{[3]}$, that is, the following holds.

\begin{corollary}\label{detremcor} \sl Every nil-polynomial $P\in\CC[W]$ is uniquely determined by its quadratic and cubic terms.
\end{corollary}

Let $P\in\CC[W]$ be a nil-polynomial. Without loss of generality we may assume that
$W=\CC^{n}$ for $n:=\dim(\5N)-1$. There exists a basis $e_1,\dots,e_n$ of $\ker\omega$ such that $\varphi(x)=\sum_{\alpha=1}^nx_{\alpha}e_{\alpha}$ for $x=(x_{1},\dots,x_{n})\in\CC^{n}$, and we write $\CC[W]=\CC[x_1,\dots,x_n]$. Then
$$
P^{[2]}(x)=\sum_{\alpha,\beta=1}^ng_{\alpha\beta}x_{\alpha}x_{\beta}\,,\qquad
P^{[3]}(x)=\sum_{\alpha,\beta,\gamma=1}^nh_{\alpha\beta\gamma}x_{\alpha}x_{\beta}x_{\gamma}\,,
$$ 
where $g_{\alpha\beta}$ and $h_{\alpha\beta\gamma}$ are symmetric in all indices. As stated above, for the admissible projection $\pi$ with $\ker\pi=\ker\omega$, the hypersurface $S_{\pi}$ is linearly equivalent to the graph
\begin{equation}
S:=\{(x_{0},x)\in\CC\times\CC^{n}:x_0=P(x)\}\,.\label{eqs}
\end{equation}

\begin{proposition}\label{Blaschke} \sl The equation of $S$ in
  {\rm(\ref{eqs})} is in Blaschke normal form, that is,
  $\sum_{\alpha\beta}g^{\alpha\beta}h_{\alpha\beta\gamma}=0$ for all
  $\gamma$, where $(g^{\alpha\beta}):=(g_{\alpha\beta})^{-1}$.
\end{proposition}

\noindent {\bf Proof:} Choose a mapping $\tau:\{1,\dots,n\}\ra\NN$
such that for every $j\ge1$ the set $\{e_{\alpha}:\tau(\alpha)=j\}$
gives a basis of ${\mathcal N}^j/{\mathcal N}^{j+1}$. Clearly,
$g_{\alpha\beta}=0$ if $\tau(\alpha)+\tau(\beta)>\nu$. This implies that
$g^{\alpha\beta}=0$ if $\tau(\alpha)+\tau(\beta)<\nu$. Since
$a_{\alpha\beta\gamma}=0$ if $\tau(\alpha)+\tau(\beta)\ge \nu$, the
proposition follows.\qed \vspace{0.3cm}

As shown in Proposition 1 of \cite{EE} (and also
\cite{L}), every linear isomorphism between hypersurfaces in Blaschke normal form is written as
$$
x_0\mapsto\hbox{const}\cdot x_0\,,\quad x\mapsto Cx\,,
$$
where $C\in\GL(n,\CC)$. Hence Propositions \ref{prop1prime}, \ref{Blaschke} and Corollaries \ref{lineqalg}, \ref{detremcor} yield the following theorem, which is the main result of this section.

\begin{theorem}\label{main1}\sl Let $P,\tilde P\in\CC[x_1,\dots,x_n]$ be arbitrary
  nil-polynomials arising from admissible algebras $\5N$,
  $\tilde{\5N}$, and let $S$, $\tilde S\subset\CC^{n+1}$ be the graphs of $P$, $\tilde P$, respectively. Then if at least one of the hypersurfaces $S$, $\tilde S$ is affinely homogeneous (e.g. if one of $\5N$, $\tilde{\5N}$ admits a grading), the following
  conditions are equivalent: \0 $\5N$, $\tilde{\5N}$ are isomorphic as associative algebras,\1$S$, $\tilde S$ are affinely equivalent, \1$S$, $\tilde
  S$ are linearly equivalent, \1$P$, $\tilde P$ are linearly equivalent, \1there exist $c\in\CC^{*}$ and $C\in\GL(n,\CC)$ with
  \begin{equation}
  c{\cdot}\tilde P^{[\ell]}(x)=P^{[\ell]}(Cx), \quad \ell=2,3 \label{identsk}
  \end{equation}
\phantom{XXX} for all $x\in\CC^{n}$.
\end{theorem}

We finish this section with an example illustrating Theorem \ref{main1} for small values of $n$.

\begin{example}\label{equivnilpolc}\rm For $n=0$ there is only one nil-polynomial, namely $0$, corresponding to the only 1-dimensional admissible algebra up to isomorphism. For $n=1$ there is again only one linear equivalence class of nil-polynomials represented by $x_1^2$; it corresponds to the only 2-dimensional admissible algebra up to isomorphism, which is generated by an element of order 3. For $n=2$ there are two linear equivalence classes of nil-polynomials represented by $P:=x_{1}x_{2}+x_{1}^{3}$ and $\tilde P:=x_1x_2$. The nil-polynomial $P$ corresponds to the isomorphism class of the 3-dimensional admissible algebra generated by an element of order 4, whereas $\tilde P$ corresponds to the isomorphism class of the 3-dimensional admissible algebra obtained from two 2-dimensional admissible algebras by identifying their annihilators. 
\end{example}

\section{Application to Quasi-Homogeneous\\ Isolated Hypersurface Singularities}\label{section2}
\setcounter{equation}{0}

In this section we apply our results of Section \ref{section1} to the biholomorphic equivalence problem for isolated hypersurface singularities. Let $V$ be a hypersurface germ at the origin in $\CC^m$ having an isolated singularity (hence $m\ge 2$), and ${\mathcal N}(V)$ the maximal ideal of the moduli algebra ${\mathcal A}(V)$ of $V$. Clearly, the unital extension ${\mathcal N}(V)^0$ of ${\mathcal N}(V)$ is isomorphic to ${\mathcal A}(V)$. By the Mather-Yau theorem, ${\mathcal N}(V)=\{0\}$ if and only if $V$ is biholomorphically equivalent to the germ of the hypersurface $z_1^2+\dots+z_m^2=0$ at the origin. Thus for our purposes it is sufficient to assume that ${\mathcal N}(V)$ is non-zero. 

We consider quasi-homogeneous singularities (see the introduction for the definition). The following is well-known, but for completeness of our exposition we include a proof.

\begin{proposition}\label{sinadmiss}\sl Let $V$ be a hypersurface germ in $\CC^m$ having an isolated singularity, and assume that ${\mathcal N}(V)$ is non-zero. Then the singularity of $V$ is quasi-homogeneous if and only if ${\mathcal N}(V)$ is a complex admissible algebra.
\end{proposition}

\noindent {\bf Proof:} Suppose first that the singularity of $V$ is quasi-ho\-mo\-ge\-neous. For every generator $f$ of $I(V)$ we have ${\mathcal A}(V)={\mathcal O}_m/J(f)$ (see the introduction). Therefore ${\mathcal A}(V)$ is a complete intersection ring, which implies that ${\mathcal N}(V)$ is an admissible algebra (see \cite{B}). Conversely, it follows from \cite{K} (see also Remark (3.7) in \cite{A}) that if the singularity of $V$ is not quasi-homogeneous, one has $\dim\Ann({\mathcal N}(V))\ge 2$, thus ${\mathcal N}(V)$ is not admissible.\phantom{XXX}\qed
\vspace{0.3cm}

\begin{remark}\label{gorenstci}\rm Not every graded admissible algebra can be realized as the maximal ideal of the moduli algebra of a quasi-homogeneous isolated hypersurface singularity. A simple example is provided by the nil-polynomial $x_{1}x_{2}+x_{3}^{2}$ on $\CC^{3}$, which arises from the graded admissible algebra $\5N=\CC^{4}$ with product $(u_{1},u_{2},u_{3},u_{4})\cdot(v_{1},v_{2},v_{3},v_{4}):=(0,0,0,u_{1}v_{2}+u_{2}v_{1}+2u_{3}v_{3})$.  The unital extension $\5N^{\,0}$ is isomorphic to ${\mathcal O}_3 /{\mathfrak I}$ with $\mathfrak I:=(z_1^2,z_2^2,z_1z_3,z_2z_3,z_1z_2+z_3^2)$. Since the minimal number of generators of the ideal $\mathfrak I$ is five, $\5N^{\,0}$ is not even a complete intersection ring.
\end{remark}

\smallskip

\begin{theorem}\label{main2} \sl Let $V$, $\tilde V$ be
  hypersurface germs in $\CC^m$ each having a quasi-homogeneous
  isolated singularity, and assume that $\5N(V)$, $\5N(\tilde V)$ are non-zero. Let furthermore $P,\tilde P\in\CC[x_{1},\dots,x_{n}]$
  be arbitrary nil-polynomials arising from the admissible algebras $\5N(V)$, $\5N(\tilde V)$,
  respectively. Then the germs $V$, $\tilde V$ are biholomorphically
  equivalent if and only if the nil-polynomials $P$, $\tilde P$ are linearly equivalent, that is, if $c{\cdot}\tilde P(x)=P(Cx)$ for all
  $x\in\CC^{n}$ and suitable $c\in\CC^{*}$, $C\in\GL(n,\CC)$. This occurs if and only if identities {\rm(\ref{identsk})} hold.
\end{theorem}

\noindent{\bf Proof:} If $V$, $\tilde V$ are biholomorphically
equivalent, ${\mathcal N}(V)$, ${\mathcal N}(\tilde V)$ are clearly
isomorphic. Conversely, if ${\mathcal N}(V)$, ${\mathcal N}(\tilde
V)$ are isomorphic, $V$, $\tilde V$ are biholomorphically
equivalent by the Mather-Yau theorem. Since the singularity of each of
$V$, $\tilde V$ is quasi-homogeneous, the algebras ${\mathcal N}(V)$, ${\mathcal N}(\tilde V)$ 
admit gradings. The theorem now follows
from Theorem \ref{main1}.\qed \vspace{0.3cm}

Theorem \ref{main2} shows that the verification of biholomorphic equivalence of
two hypersurface germs $V$, $\tilde V$ in $\CC^m$ having
quasi-homogeneous isolated singularities reduces to the
invariant theory for pairs of quadratic and cubic forms in at most
$\dim{\5A}(V){-}2$ variables. On the other hand, for the verification of
biholomorphic non-equivalence of $V$, $\tilde V$ it is especially
useful to exploit higher-order terms in identities (\ref{identsk}) since they depend on fewer variables.

In the following example we demonstrate how Theorem \ref{main2} can be used
to distinguish biholomorphically non-equivalent simple elliptic hypersurface
singularities. Recall that all such singularities split into three types $\tilde E_6$, $\tilde E_7$, $\tilde E_8$ (see \cite{Sa2}). Biholomorphic equivalence of simple elliptic singularities is well-understood (see \cite{Sa2}, \cite{CSY}, \cite{SY}, \cite{Ea}), and therefore we restrict ourselves to the most interesting case of $\tilde E_8$-singularities.

\begin{example}\label{exe8}\rm Elliptic singularities of type $\tilde E_8$ are the quasi-homogeneous singularities at the origin of the following hypersurfaces in $\CC^3$:
$$
V_t:=\left\{(z_1,z_2,z_3)\in\CC^3:z_1^6+tz_1^4z_2+z_2^3+z_3^2=0\right\},
$$
where $t\in\CC$ satisfies $4t^3+27\ne 0$. Then the statement
\begin{equation}
\hbox{the germs of $V_{r}$ and $V_s$ at 0 are biholomorphic}\quad\Longleftrightarrow\quad r^{3}=s^{3}\label{relsss}
\end{equation}
is well-known. The transformation
$(z_{1},z_{2},z_{3})\mapsto(z_{1},\rho z_{2},z_{3})$, with $\rho^3=1$, maps $V_{t}$
biholomorphically to $V_{\rho t}$, thus the implication
$\Longleftarrow$ in (\ref{relsss}) is trivial. An elementary
proof of the converse implication is as follows.

As in \cite{SY}, \cite{Ea}, we consider the monomials
$$
z_1^4z_2,\,\, z_1,\,\, z_2,\,\, z_1^2,\,\,z_1z_2,\,\,z_1^3,\,\,z_1^2z_2,\,\,z_1^4,\,\,z_1^3z_2,
$$
and let $e_0,\dots,e_8$ be the vectors in ${\mathcal N}_t:={\mathcal N}(V_t)$ arising from them. These vectors form a basis of ${\mathcal N}_t$. Let $\pi_t$ be the projection on ${\mathcal N}_t$ with range $\Ann({\mathcal N}_t)=\langle e_0\rangle$ such that $\ker\pi_t$ is spanned by $e_1,\dots,e_8$. From these data the following nil-polynomial in $\CC[x_{1},\dots,x_{8}]$ is derived:
$$
\begin{array}{l}
\displaystyle P_t:=-\frac{t}{1080}x_1^6+\frac{1}{72}x_1^4\left(3x_2-2tx_3\right)\;+\\
\vspace{-0.1cm}\\
\hspace{1cm}\displaystyle\frac{1}{18}x_1^2\left(3x_1x_4-2tx_1x_5+t^2x_2^2+9x_2x_3-3tx_3^2\right)\;-\\
\vspace{-0.1cm}\\
\hspace{1cm}\displaystyle\frac{1}{18}Q_t+\hbox{cubic terms involving $x_1$}+\hbox{quadratic terms}\,,\\
\end{array}
$$
where
\begin{equation}
Q_t:=tx_2^3-2t^2x_2^2x_3-9x_2x_3^2+2tx_3^3\;.\label{qt}
\end{equation}

Suppose that for some $r\ne s$ the germs of $V_{r}$ and $V_{s}$ are biholomorphically equivalent. Since 0 is the only value of $t$ for which $P_t$ has degree 6, we have $r,s\ne 0$. By Theorem \ref{main2} there exist $c\in\CC^*$ and $C\in\GL(8,\CC)$ such that   
\begin{equation}
c\cdot P_r(x)\equiv P_s(Cx)\,.\label{idene8}
\end{equation}
By comparing terms of order 6 in identity (\ref{idene8}), we obtain that the first row in the matrix $C$ has the form $(\mu,0,\dots,0)$, and
\begin{equation}
c=\frac{s}{r}\mu^6.\label{eq1}
\end{equation}
Next, let $({}*{},\alpha,\beta,{}*{},\dots,{}*{})$ and $({}*{},\gamma,\delta,{}*{},\dots,{}*{})$ be the second and third rows in $C$, respectively, for some $\alpha,\beta,\gamma,\delta\in\CC$. Comparing the terms of order 4 in (\ref{idene8}) that do not involve $x_1^3$, we see that the matrix
$$
D:=\left(
\begin{array}{ll}
\alpha & \beta\\
\gamma & \delta
\end{array}
\right)
$$
is non-degenerate. Further, comparing terms of order 5 in (\ref{idene8}) we obtain
\begin{equation}
\displaystyle\beta=\frac{2}{9}\left(-3\alpha r+3\delta s+2\gamma rs\right)\label{beta}
\end{equation} 
and
\begin{equation}
\displaystyle c=\left(\alpha-\frac{2s}{3}\gamma\right)\mu^4.\label{eq2}
\end{equation}

We will now compare the terms of order 3 in (\ref{idene8}) that depend only on $x_2$, $x_3$. For $X:=(x_2,x_3)$
we have
\begin{equation}
c\cdot Q_{r}(X)=Q_{s}(DX).\label{qttt}
\end{equation}
Setting 
$$
\displaystyle D_t:=\left(
\begin{array}{cc}
1/3 & 2t/3\\
\vspace{-0.3cm}\\
0 & 1
\end{array}
\right),
$$
one observes
$$
\displaystyle Q_t(D_tX)={\mathbf Q}_t(X):=\frac{t}{27}x_2^3-3\Delta_tx_2x_3^2-4t\Delta_tx_3^3\;,
$$
where $\Delta_t:=1+4t^3/27$. Hence (\ref{qttt}) implies
\begin{equation}
c\cdot{\mathbf Q}_{r}(X)={\mathbf Q}_{s}(\hat DX)\,,\label{tildeq}
\end{equation}
where $\hat D:=D_{s}^{-1}D D_{r}$. By (\ref{beta}) we have
$$
\hat D=\left(
\begin{array}{ll}
a & 0\\
b & d
\end{array}
\right),
$$
with $a:=\alpha-2s\gamma/3$, $b:=\gamma/3$, $d:=\delta+2r\gamma/3$.

It follows from (\ref{tildeq}) and the non-degeneracy of $\hat D$ that $b(a+2sb)=0$. If $b=0$, comparison of the three pairs of coefficients in (\ref{tildeq}) yields
$$
\displaystyle c=\frac{s}{r}a^3=\frac{\Delta_{s}}{\Delta_{r}}ad^2=\frac{s\Delta_{s}}{r\Delta_{r}}d^3\,.
$$
Therefore $r^3\Delta_{s}=s^3\Delta_{r}$, and we obtain that in this case $r^3=s^3$. Let now $b\ne 0$, that is, $a=-2sb$. In this situation comparison of the three pairs of coefficients in (\ref{tildeq}) yields
\begin{equation}
\displaystyle c=54\frac{s}{r}b^3=2\frac{s\Delta_{s}}{\Delta_{r}}bd^2=\frac{s\Delta_{s}}{r\Delta_{r}}d^3.\label{eq3}
\end{equation}
From identities (\ref{eq1}), (\ref{eq2}) and the first equality in (\ref{eq3}) we obtain $\Delta_{r}=0$, contrary to the initial assumption on the parameter $t$. This shows that $r^3=s^3$, as required.\qed
\end{example}

We mention that the family of cubic polynomials $Q_t$ defined in (\ref{qt}) was introduced in \cite{Ea}, where a proof of (\ref{relsss})  relying on the invariant theory for $Q_t$ with respect to upper-triangular matrices was given. In our approach we do not need to refer to any invariant theory, it is sufficient to perform elementary manipulations with certain homogeneous polynomials of orders three to six.
\vspace{0.3cm}

We now give another example of application of Theorem \ref{main2}.

\begin{example}\label{newexample}\rm Consider the following family of curves in $\CC^2$:
$$
V_t:=\left\{(z_1,z_2)\in\CC^2:z_1^4+tz_1^2z_2^3+z_2^6=0\right\},
$$
where $t\in\CC$. The curve $V_t$ has a quasi-homogeneous isolated singularity at the origin if and only if $t^2\ne 4$. We show:
\begin{equation}
\hbox{the germs of $V_{r}$, $V_s$ at 0 are biholomorphic}\quad\Longleftrightarrow\quad r^{2}=s^{2}.\label{relssss}
\end{equation}
The implication $\Longleftarrow$ in (\ref{relssss}) is trivial since the map $(z_{1},z_{2})\mapsto(z_{1},-z_{2})$ transforms $V_t$ into $V_{-t}$.

For the converse implication we assume that $t^2\ne 4$ and consider the following monomials:
$$
z_1^2z_2^4,\,\, z_2,\,\,z_1,\,\,z_1^2,\,\,z_1z_2,\,\,z_2^2,\,\,z_1^2z_2,\,\,z_1z_2^2,\,\,z_2^3,\,\,z_1z_2^3,\,\,z_1^2z_2^2,\,\,z_2^4,\,\,z_1^2z_2^3,\,\,z_1z_2^4\,.
$$
Let $e_0,\dots,e_{13}$ be the vectors in ${\mathcal N}_t:={\mathcal N}(V_t)$ arising from them. These vectors form a basis of ${\mathcal N}_t$. Let $\pi_t$ be the projection on ${\mathcal N}_t$ with range $\Ann({\mathcal N}_t)=\langle e_0\rangle$ such that $\ker\pi_t$ is spanned by $e_1,\dots,e_{13}$. From these data the following nil-polynomial in $\CC[x_1,\dots,x_{13}]$ is derived:
$$
\begin{array}{l}
\displaystyle P_t:=-\frac{t}{10080}x_1^7+\frac{1}{48}x_1^4\left(x_2^2-\frac{t}{5}x_1x_5\right)-\frac{t}{48}x_1x_2^4+\frac{1}{4}x_1^2x_2^2x_5+\frac{1}{6}x_1^3x_2x_4\;-\\
\vspace{-0.3cm}\\
\displaystyle\hspace{1cm}\frac{t}{24}x_1^3x_5^2-\frac{t}{48}x_1^4x_8+\frac{1}{24}x_1^4x_3+\hbox{terms of degree $\le 4\,$}.
\end{array}
$$

Suppose that for some $r\ne s$ the germs of $V_{r}$ and $V_{s}$ are biholomorphically equivalent. Since 0 is the only value of $t$ for which $P_t$ has degree 7, we have $r,s\ne 0$. By Theorem \ref{main2} there exist $c\in\CC^*$ and $C\in\GL(13,\CC)$ such that   
\begin{equation}
c\cdot P_r(x)\equiv P_s(Cx)\,.\label{equivnew}
\end{equation}
Comparing terms of order 7 in identity (\ref{equivnew}), we obtain that the first row in the matrix $C$ has the form $(\mu,0,\dots,0)$, and
\begin{equation}
c=\frac{s}{r}\mu^7\,.\label{7}
\end{equation}
Next, comparing terms of order 6 in (\ref{equivnew}), we see that the second row in the matrix $C$ has the form $(\alpha,\beta,0,\dots,0)$,  and 
\begin{equation}
c=\mu^4\beta^2\,.\label{6}
\end{equation}
Further, comparing the terms of order 5 in (\ref{equivnew}) that do not involve $x_1^2$, we obtain
\begin{equation}
c=\frac{s}{r}\mu\beta^4\,.\label{5}
\end{equation}
From (\ref{7}), (\ref{6}), (\ref{5}) we get $r^2=s^2$, as required. \qed

\end{example}

{\obeylines
G. Fels and W. Kaup:
Mathematisches Institut
Universit\"at T\"ubingen
Auf der Morgenstelle 10
72076 T\"ubingen
Germany
e-mail: {\tt gfels@uni-tuebingen.de, kaup@uni-tuebingen.de}
\hbox{ \ \ }
A. Isaev:
Department of Mathematics
The Australian National University
Canberra, ACT 0200
Australia
e-mail: {\tt alexander.isaev@anu.edu.au}
\hbox{ \ \ }
N. Kruzhilin:
Department of Complex Analysis
Steklov Mathematical Institute
8 Gubkina St.
Moscow GSP-1 119991
Russia
e-mail: {\tt kruzhil@mi.ras.ru}
}

\end{document}